%% file: frame-rokhlin.tex
\documentclass{amsart}
\usepackage{amssymb,graphicx,rotating,multirow,multicol}
\usepackage{amsmath,amsthm,amsfonts,amscd,accents}
\usepackage[toc,page]{appendix}
\usepackage{caption,comment,slashed}
\usepackage[mathscr]{euscript}
\usepackage{epstopdf}
\usepackage{float}
\usepackage{epigraph}
\usepackage{epsf}
\usepackage{geometry}
\usepackage[symbol]{footmisc}

\input xy
\xyoption{all}

\swapnumbers
\newtheorem{theorem}{Theorem}[subsection]
\newtheorem{lemma}[theorem]{Lemma}

\newtheorem{cor}[theorem]{Corollary}
\newtheorem{proposition}[theorem]{Proposition}
\theoremstyle{remark}

\numberwithin{equation}{subsection}

{\catcode`\@=11 \gdef\mnote#1{\marginpar{\tiny
 \tolerance\@M\spaceskip2.6\p@ plus10\p@ minus.9\p@\rm#1}}}

\let\Bbb\mathbb

\def\s{\sigma}
\def\sm{\smallsetminus}

\newcommand{\be}{\begin{equation}}
\newcommand{\ee}{\end{equation}}
\let\ge\geqslant 
\let\le\leqslant 

\let\til\widetilde

\let\a\alpha

\def\Z{\Bbb Z}
\def\R{\Bbb R}
\def\C{\Bbb C}
\def\H{\Bbb H}

\def\S{\Sigma}

\def\Cp#1{\operatorname{\Bbb CP^{#1}}}
\def\barCp#1{\operatorname{\overline{\Bbb CP}^{#1}}}
\def\Rp#1{\Bbb{RP}^{#1}}

\def\SO{\operatorname{SO}}

\def\Spin{\operatorname{Spin}}

\def\SU{\operatorname{SU}}

\let\+\sqcup
\let\-\smallsetminus
\def\ind{\operatorname{ind}}
\def\Ind{\operatorname{Ind}}

\def\pt{\operatorname{pt}}
\def\M{\accentset{\circ}M}
\def\X{\accentset{\circ}X}

\def\NN{L}

\def\card{\operatorname{card}}

\makeatletter
\newcommand{\addresseshere}{%
  \enddoc@text\let\enddoc@text\relax
}
\makeatother

\begin{document}


\title[A glimpse into Rokhlin's Signature Divisibility Theorem]
{A glimpse into Rokhlin's Signature Divisibility Theorem}

\author[]
{S.~Finashin, V.~Kharlamov}

\begin{abstract}
In the main section
we give an overview of  Rokhlin's proof of his famous
theorem on divisibility of signature by 16. In the appendix we retrace
some of further developments that show how this theorem became a cornerstone
in the contemporary theory of manifolds.
\end{abstract}

\maketitle

\setlength\epigraphwidth{.5\textwidth}
\epigraph{
Die reine Mathematik w\"achst, indem man alte Probleme mit neuen  Methoden durchdenkt.\\
 In dem Maasse, wie wir die fr\"uheren Aufgaben besser verstehen, bieten sich neue von selbst.}
{Felix Klein, Riemann und seine Bedeutung f\"ur die Entwickelung der modernen Mathematik.\footnotemark[2]}
%
\footnotetext[2]{\tiny``Pure mathematics advances as one thinks on old problems with new methods. As we understand the earlier tasks better, new ones come naturally.''
/ Felix Klein, {\sl Riemann and its impact on the development of modern mathematics.}}
 \renewcommand*{\thefootnote}{\arabic{footnote}}

%

%
%
\input intro.tex

\input proof.tex

\input add.tex

\end{document}

%% file: intro.tex
\section{Instead of Introduction}
\subsection{Rokhlin's Signature Divisibility Theorem}
This theorem can be stated as follows.

\begin{theorem}
[\cite{R4}]\label{divisibility}
The signature of any closed oriented smooth $4$-manifold $M$ with $w_2=0$ is divisible by $16$.
\end{theorem}

It was published in the last one of a series of four
short notes, \cite{R1}--\cite{R4}, that appeared in 1951--52.
The initial motivation was calculation of the homotopy groups $\pi_{n+3}(S^n)$.
But the output went far beyond that
narrow
scope: intrinsic (versus Pontryagin's framed) cobordism groups have emerged, fine geometric methods for their calculation were developed and then applied for calculation of cobordism groups in dimensions 3 and 4, and as a by-product the Signature Formula (Theorem \ref{formula}) and the Signature Divisibility Theorem (Theorem \ref{divisibility}) were established.

The fate of Rokhlin's proofs of his theorems
on four-dimensional manifolds, and first and foremost the theorem on divisibility of the signature by 16, was unusual. While the results themselves were acknowledged at once and became widely used,
the original proofs of Rokhlin were
apparently
not understood for a long time.
At first glance this happened because the published text of the proofs was very concise,
and their understanding was complicated by an error made in the first of the four notes and corrected in the fourth. But the main reason is deeper.
Around the same time, i.e., at the beginning of the 1950's, a new powerful method for the study of homotopy groups ("the method of killing homotopy groups") was elaborated by H.~Cartan and J.-P.~Serre and successfully applied by  J.-P.~Serre \cite{SerreNote}
to calculation of the homotopy groups $\pi_{n+3}(S^n)$. This opened a way for applying in the reverse direction
the connection between homotopy problems and geometric problems that was
 discovered by Pontryagin, and exploited by Pontryagin and Rokhlin for solving homotopy problems.
In particular, it became easier to derive Rokhlin's theorem from computation of  $\pi_{n+3}(S^n)$,
although originally this theorem was obtained by Rokhlin as a by-product (if not a lemma) in the process of computation of these groups. Furthemore, the new methods of the homotopy theory allowed R.~Thom \cite{ThomPaper} to go much further in understanding of the cobordism groups in all dimensions. By all these reasons Rokhlin abandoned publication of a detailed text with the original proofs.

Since Serre's computation of $\pi_{n+3}(S^n)$ 
looked more elementary,
it became natural to deduce Rokhlin's divisibility theorem from that computation.
It
was done by M.~Kervaire and J.~Milnor in \cite{KM2}
and then was reproduced in many other publications in 50th--60th. This resulted in
a widely-spread misconception on the nature of Rokhlin's original proof
(Rokhlin himself gave in \cite{Kolomna} an additional impetus to this confusion by offering a "short proof" of the divisibility theorem, also deducing it from Serre's homotopy calculation of  $\pi_{n+3}(S^n)$).


It is only in the 1970's that an understanding of the ubiquitous
role of Rokhlin's Theorem \ref{divisibility}
and the intrinsic demands of four-dimensional topology led first to a reinvention of the original geometric proof, and then to a true reading of Rokhlin's papers
(see the comments by L. Guillou and A. Marin on Rokhlin's four notes in their book \cite{Book}).


Nevertheless
even in articles of the most respectable authors the confusion about Rokhlin's initial proof is still reproduced again and again,
so, we decided to apply extra efforts to clarify this story.
That is why we make an attempt to reconstruct the original Rokhlin's proof in its
essential details (where we mostly follow L.~Guillou and A.~Marin \cite{Book}),
and demonstrate that this proof
consisted not in
deducing Theorem \ref{divisibility} from calculation of the homotopy groups $\pi_{n+3}(S^n)$,
but
in giving both simultaneously (cf., Section \ref{conclusion}).

\subsection{Reconstruction of Rokhlin's proofs}
We tried to summarize the principal ideas from the
four Rokhlin's notes and clarify some interesting details omitted in the original as ``obvious enough''.
We made preference to the language of framings and obstructions used by Rokhlin, although from modern perspective,
it would be more simple to use language of Spin structures, classifying spaces, etc. (techniques which appeared after Rokhlin's work).

We restricted ourselves only with proving Theorem \ref{divisibility} and just stated two other fundamental Rokhlin's Theorems
\ref{cobordism} and \ref{formula} used in the proof.

In Appendix we point a few results from \cite{R1} - \cite{R4} representing an additional output of Rokhlin's proof of Theorem \ref{divisibility},
discuss further refinements and generalizations of this theorem, and
indicate some of their most striking applications
which demonstrate how essential and ubiquitous this Rokhlin's theorem has become for the topology of manifolds.

\vskip0.05in
{\it Acknowledgements.}
This text has resulted from our work on a short essay that we wrote for a book about Petersburg mathematicians (conceived on
the occasion of ICM-2022). We thank I.~Kalinin who involved us into this book project.

%% file: proof.tex
\section{Outline of Rokhlin's proof of divisibility of the signature by 16}\label{onlineversion}
\setlength\epigraphwidth{.6\textwidth}
\epigraph{"The theory of cobordism was initiated by
L.~Pontryagin and V.A.~Rokhlin. It came of age with the work of R. Thom."}
{John Milnor, {\it A survey of cobordism theory.}}

\subsection{Two basic theorems}
Rokhlin's initial proof of Theorem \ref{divisibility} was closely tied with
his calculation of the oriented cobordism groups $\Omega_3$ and $\Omega_4$, and the techniques he developed for that.

\begin{theorem}[Rokhlin's Cobordism  Groups Theorem, \cite{R3, R4}]\label{cobordism}
\quad
\begin{enumerate}
\item $\Omega_3=0${\rm:}
any closed oriented smooth 3-manifold bounds a
compact oriented smooth 4-manifold.
\item $\Omega_4=\Z${\rm:}
two closed oriented smooth $4$-manifolds are cobordant if and only if they have the same signature.
\qed\end{enumerate}
\end{theorem}

The next theorem is an immediate corollary of Theorem \ref{cobordism} (2).

\begin{theorem}[Rokhlin's Signature Formula, \cite{R4}]\label{formula}
For any closed oriented smooth $4$-manifold, the Pontryagin number $p_1[M]$  is equal to
$3\,\sigma(M)$.
\end{theorem}
\begin{proof}
Due to cobordism invariance of characteristic numbers (see \cite{P2}), $p_1[M]$ is a cobordism invariant.
By the Poincar\'e-Lefschetz duality, $\sigma(M)$ is also a cobordism invariant. Both are additive under disjoint union,
hence, it remains to notice
that $\sigma(\C P^2)=1$, $p_1[\C P^2]=3$, and to apply then Theorem \ref{cobordism}.
\end{proof}

\subsection{Framed spheres}\label{framed-spheres}
Apart of a systematic use of the geometric tools developed for proving Theorem \ref{cobordism},
Rokhlin's proof of Theorem \ref{divisibility} was based on Pontryagin's construction \cite{P1}, which establishes a canonical {\it Pontryagin isomorphism}
$\pi_{n+k}(S^n)\cong \Pi^k_n$ between the homotopy groups of spheres and groups of
normally framed cobordisms.

Recall that for an $m$-dimensional oriented manifold $M$ smoothly embedded into $\R^{m+n}$ (or equivalently, to $S^{n+m}$)
its {\it normal $n$-framing} is a trivialization of the normal bundle
of $M$ in $\R^{m+n}$
(if this bundle is trivializable)
given by $n$ pointwise linearly independent sections
and compatible with the orientations of $M$ and $\R^{m+n}$.
A {\it stabilization} of such a framing is an $(n+1)$-framing obtained by
composing embedding of $M\subset\R^{m+n}$ with $\R^{m+n}\subset\R^{m+n+1}$ and adding the constant vector field $e_{n+1}\in \R^{m+n+1}$
to the $n$-framing.
 For $n\ge m+2$ the normal $n$-framings are {\it stable}: stabilization does not change the relation of homotopy equivalence of framings.
A normal framing is said to be {\it matching} with a tangent framing
respecting the orientation of $M$ if their union
is stable homotopic to the
standard (constant)
tangent framing of $\R^{m+n}$ restricted to $M$.
The matching relation gives a one-to-one correspondence after we pass to the stable homotopy classes
of normal and tangent framings.

For the standard sphere $S^m\subset\R^{m+1}\subset\R^{m+n}$,
the normal bundle of $S^m$
has a canonical
normal n-framing $(\alpha_1,\alpha_2,\dots , \alpha_n)$ where $\alpha_1$ is the outward normal vector field and $\alpha_{k+2}$ with $k\ge 0$ are the coordinate vectors $e_{m+k+2}$,
so that according to standard matrix convention
any other normal $n$-framing $(\beta_1, \dots, \beta_n)$ of this embedding
can be identified with a map $F: S^m\to\SO_{n}$, $x\in S^m\mapsto \{f_{ij}(x)\}\in\SO_n$ as follows:
\begin{equation}\label{functor}
\beta_j(x)=\sum \alpha_i(x) f_{ij}(x)  \quad (\beta = \alpha F   \text{  in short}).
\end{equation}
Thus, Pontryagin's construction defines a homomorphism
\footnote{We switched from Rokhlin's notation $D$ for this homomorphism to $J$, 
since as a map $\pi_m(\SO_n)\to\pi_{m+n}(S^n)$ it is the standard $J$-homomorphism.
}
\ \ $\pi_m(\SO_n)\xrightarrow{J} \Pi^m_n\cong\pi_{m+n}(S^n)$
commuting with the
stabilization induced by the inclusion \ $\operatorname{in}:\SO_n\to\SO_{n+1}$ and the suspensions, $\Sigma$ and ${\Sigma}^{\Pi}$, of
Freudental and Pontryagin, respectively:
$$
\xymatrix{
\pi_m(\SO_n)\ar[d]^{in_*}\ar[r]^{\ J} &\Pi^m_n\ar[d]^{\Sigma^\Pi} \ar@{*{\quad}*{\quad\cong}*{\quad\quad}}[r] &\quad\pi_{m+n}(S^n)\ar[d]^\Sigma\\
\pi_m(\SO_{n+1})\ar[r]^{\ \  J}        &\hskip3mm\Pi^m_{n+1}                \ar@{*{{\quad}\quad \cong\qquad\qquad}}[r] &\,\,\pi_{m+n+1}(S^{n+1})
}
$$

A sphere $S^3\subset\H$ as a Lie group of unit quaternions has two natural
tangent framings: the left-invariant framing
$\xi(x)=(x\mathrm{i},x\mathrm{j},x\mathrm{k})$ and the right-invariant framing $\xi'(x)=(\mathrm{i}x,\mathrm{j}x,\mathrm{k}x)$.
The relation $\xi'(x)=-\overline{\xi(\bar x)}$ shows that $\xi'$ is related to $\xi$ by the quaternion conjugation.

Similar formulas define normal 4-framings $\hat\xi_4,\hat\xi'_4$ as sections of the normal bundle $S^3\times\H\to S^3$ of $S^3\subset \R^7$.
Here, the normal bundle of $S^3\subset \R^7$
is trivialized by
the consecutive
embeddings $S^3\subset \R^4\subset \R^5\subset \R^6\subset \R^7$ whose outward normal vector fields are considered as
$1,\mathrm{i},\mathrm{j},\mathrm{k}\in\H$. Then, considering as above the points of $S^3$ as unit quaternions,
we put $\hat\xi_4(x)=(x,\mathrm{i}x,\mathrm{j}x,\mathrm{k}x)$ and $\hat\xi'_4(x)=(x,x\mathrm{i},x\mathrm{j},x\mathrm{k})$.
For $n\ge 5$, by taking $(n-4)$-step stabilization of $\hat\xi_4$ and $\hat\xi'_4$ we obtain
normal $n$-framings on $S^3$, which we denote
by $\hat\xi_n$ and $\hat\xi'_n$.
The framings $\hat\xi_n, \hat\xi'_n$ with $n\ge 4$ determine some elements $g_n, g_n'$ in the Pontryagin group $\Pi^3_n$ and provide also elements $b_n, b'_n\in\pi_3(\SO_n)$ such that  $J(b_n)=g_n$ and $J(b_n')=g_n'$.

Consider also a generator $a_3\in \pi_3(\SO_3)$ representing the double covering $S^3\to\Rp3=\SO_3$, which
sends $x\in S^3\subset\H$ to the rotation $h\mapsto \phi_x(h)=xhx^{-1}$, where
$h\in \R^3\subset\H$ is an imaginary quaternion.
We let $a_n\in\SO_n$, $n\ge4$, be the $(n-3)$-step stabilization of $a_3$
put $f_n=J(a_n)\in \Pi_n^3$.

The following statement is immediate from definitions.
\begin{lemma}\label{spherical-classes}
For $n\ge5$, $b_n=-b'_n$ is a generator of $\pi_3(\SO_n)=\Z$ and $g_n=-g'_n$ is a generator of $J(\pi_3(\SO(n)))$, while
$a_n=2b_n$ and $f_n=2g_n$.
\qed
\end{lemma}

\subsection{Framed tori}\label{framed-tori}
On a {\it trivialized $m$-torus} $T^m=S^1\times \dots\times S^1$ with angle coordinates $\theta=(\theta_1,\dots,\theta_m)$
any Lie-group-invariant tangent framing (compatible with the orientation of $T^m$) is isotopic to
$\eta^m=(\frac\partial{\partial\theta_1},\dots ,\frac\partial{\partial\theta_m})$.
 Moreover, on a smooth oriented $m$-torus, $T$,
the homotopy class of  a Lie-group-invariant tangent framing is independent of a particular orientation-preserving trivialization $T\cong T^m$ (or introducing a Lie-group structure in any other way).

We define a {\it standard embedded torus} $T^m\subset \R^{m+1}\subset S^{m+1}$ inductively:
starting from $T^1=S^1\subset\R^2$,
we define then $T^m$ as the boundary of a tubular neighborhood $N(T^{m-1})\cong T^{m-1}\times D^2$ of the composed embedding $T^{m-1}\subset\R^m\subset\R^{m+1}$.
Let $\nu^m_1$ denote the unit outward-normal vector field along $T^m$ and $t=(t_1,\dots,t_m)$ the orthonormal tangent vector field ($t_i$ is obtained by
scaling of $\frac\partial{\partial\theta_i}$ that makes it a unit vector).

Let $r(\phi)\in\SO_2$ denote the matrix of rotation by angle $\phi\in S^1$, and let $r_{i,j}(\phi)\in\SO_n$, $1\le i<j\le n$, denote
the image of $r(\phi)$ under the homomorphism
$\SO_2\to\SO_N$ induced by the inclusion of the coordinate $ij$-plane $\R^2\subset\R^n$.
Then, in accordance with convention (\ref{functor}) for matrix frame action, the framing $f=(\nu^m_1,t_1,\dots,t_m)$ is related to the (constant) Cartesian framing $e=(e_1,\dots,e_{m+1})$ in $\R^{m+1}$ by
a unique map $R : T^m \to \SO_{m+1}$ that gives $f(\theta)=e R(\theta)$.
The latter action can be decomposed
in a sequence of plane rotations defined by the following matrix product factorization
\begin{equation}\label{rotation}
R(\theta)= r_{1,2}(\theta_1)r_{1,3}(\theta_{2})\dots r_{1,m+1}(\theta_m).
\end{equation}

For a composite embedding $T^m\subset\R^{m+1}\subset\R^{m+n}$ with $n\ge 2$, we
define a {\it twisted normal framing}

\begin{equation}\label{h-framing}
\hat\eta^m_n(\theta)=
(\nu^m_1,e_{m+2}, \dots, e_{m+n}) r_{1,2}(-\theta_1-\dots-\theta_m)
\end{equation}
Each of the framings $\hat\eta_n^m$ defines a normal framed cobordism class
$h_n^m\in\Pi^m_n$ and the corresponding by the Pontryagin isomorphism homotopy class
$[h_n^m]\in\pi_{n+m}(S^n)$. By definition,
\begin{equation}\label{suspension}
h_{n+1}^{m}=\Sigma^\Pi(h_n^m)\quad \text {and} \quad  [h_{n+1}^m]=\Sigma([h_n^m]).
\end{equation}

Another way to describe the elements $h_n^m$, $m\ge1$, $n\ge2$, is to consider a framed torus $T^m\subset S^{m+n}$
given by Pontryagin's construction applied
to the composition
\begin{equation}\label{composition}
S^{m+n}\xrightarrow{\S^{m+n-3}H}S^{m+n-1}\xrightarrow{\S^{m+n-2}H}\dots\xrightarrow{\S^{n-2}H}S^{n}
\end{equation}
of $m$ consecutive suspensions over the Hopf map $H: S^3\to S^2$. In particular, $[h_2^1]$ is the generator of $\pi_3(S^2)=\Z$ given by the Hopf fibration, $[h^1_3]=\Sigma [h_2^1]$ is the generator of $\pi_4(S^3)=\Z/2$, and
\begin{equation}
[h_n^m]=\S^{n-2} [h_2^1]\circ \S^{n-1} [h_2^1]\circ\dots\circ\S^{m+n-3} [h_2^1] \quad \text{ for any $m\ge1$, $n\ge 2$}.
\end{equation}

\begin{lemma}\label{tori-framings}
\quad
\begin{enumerate}
\item
The framings $\eta^m$ and $\hat\eta_n^m$ are matching for any $n\ge2$, $m\ge1$.
\item $2h_n^m=0\in\Pi_n^m$ and $2[h_n^m]=0\in\pi_{n+m}(S^n)$ for any $n\ge 3$, $m\ge 1$.
\end{enumerate}
\end{lemma}

\begin{proof}
To prove (1), it is sufficient to check the matching for $n=2$. Due to (\ref{rotation}) and (\ref{h-framing}), a
combined framing $f^{comb}=(n_1^m,t_1,\dots,t_m,n_2^m)$, where $(n_1^m(\theta),n_2^m(\theta))=\hat\eta^m_2(\theta)$,
is obtained from the constant framing $e=(e_1,\dots,e_{m+2})$
as
\[f^{comb}(\theta)=
e R(\theta) r_{1,m+2}(-\theta_1-\dots-\theta_m).\]

The inclusion homomorphisms $r_{i,j}:\SO_2\to\SO_{m+2}$, $m\ge1$, realize the same (non-trivial) element of $\pi_1(\SO_{m+2})=\Z/2$ and so,
all $r_{i,j}$ are homotopic to each other. It implies that the map $R :T^m\to SO_{m+2}$ is homotopic to the map $T^m\to SO_{m+2}$ that sends $\theta$ to
 $$
r_{1,m+2}(\theta_1)r_{1,m+2}(\theta_{2})\dots r_{1,m+2}(\theta_m)=r_{1,m+2}(\theta_1+\dots+\theta_m),
$$
and hence the map $\theta\mapsto R(\theta) r_{1,m+2}(-\theta_1-\dots-\theta_m))$
is null-homotopic.

To prove (2), it is sufficient to check that  $2[h_n^m]=0$. The latter vanishing property is immediate from (\ref{composition}) and the
relation $2[h_3^1]=0$, which is due to $\pi_4(S^3)=\Z/2$.
\end{proof}

In what follows, we work only with the case $m=3$ and, therefore,
skipping the upper index write simply
$\eta$ instead of $\eta^3$ and $h_{n}\in\Pi^3_n$ instead of $h_n^3$.

\subsection{Relative characteristic numbers}
Since characteristic classes of smooth manifolds can be interpreted as obstructions to certain framings, one can define a {\it relative characteristic class}
as an obstruction to extend a fixed framing along a submanifold $L\subset M$
to the whole $M$.

We will be concerned primarily with the Stiefel-Whitney classes $w_2$ and the Pontryagin classes $p_1$ in the case of
$\NN=\partial M$ for compact orientable n-manifolds $M$. For instance, the relative class
$w_2(M,\zeta)\in H^2(M,\NN;\Z/2)$ is the first obstruction to extend to $M$ a framing $\zeta$
in the tangent bundle $T_*M|_\NN$.
 If a framing $\zeta$ is given
in $T_*\NN$, we use the same notation $w_2(M,\zeta)$
meaning that $\zeta$ is extended to a framing of $T_*M|_\NN$
by an outward-normal vector field to $\NN$.
It may happen also that $\NN$ is a corner in $M$: a model example is $M=M_1\times M_2$ and $\NN=\partial M_1\times\partial M_2$.
Then a framing $\zeta$ on $T_*\NN$ can be extended to $T_*M|_{\NN}$ via two outward-normal vector fields: in $M_1\times \partial M_2$ and
in $\partial M_1\times M_2$.

It is straightforward to check that
\be\label{torus-w2}
w_2(T^2\times D^2,\eta)\in H^2(T^2\times D^2,T^2\times S^1;\Z/2)=\Z/2 \quad \text{  is the generator}.
\ee
In the other words, the core torus $T^2\times0$ is a {\it characteristic surface} which means a surface Poincare-Lefschetz dual to $w_2$.
Immediately from the definition we deduce also {\it additivity} of relative characteristic classes and numbers.
Namely, for the union $M\cup M'$ of two $n$-manifolds with common boundary $\NN$
and a framing $\zeta$ over  $\NN$
the homomorphism
$H^i(M,\NN; \Z/2)\oplus H^i(M',\bar\NN; \Z/2)\to H^i(M\cup M';\Z/2)$
induces splittings
\be\label{additivity-w}
w_i(M,\zeta)\oplus w_i(M',\zeta)\mapsto w_i(M\cup M'),
\quad
f(w)[M,\NN]+f(w)[M',\NN]=f(w)[M\cup M'],
\ee
where $f(w)$ stands for any polynomial in $w_i$ and $f(w)[\dots]$ denotes the corresponding characteristic number,
obtained by evaluation of $f(w)$ on the fundamental class $[\dots]$.

Note also that for any orientable compact 4-manifold with boundary $\NN=\partial M$ the usual Wu relation $w_2(M)\circ x= x^2$, which holds for
any $x\in H^2(M,L;\Z/2)$, implies immediately that
also
\be\label{Wu-relation}
w_2(M,\zeta)\circ x=x^2
\ee
for any tangent framing $\zeta$ of the tangent bundle of $\NN$.

The Pontryagin class $p_1(V)$  is interpreted as the first obstruction to construct $n$ sections on an n-dimensional vector bundle $V$ whose span at every point has dimension $\ge n-1$ (the corresponding space, $A_n$,
of $n\times n$-matrices of rank $\ge n-1$ is
$2$-connected and has $\pi_3(A_n)=\Z$, see \cite{R5}).
If for a compact 4-manifold $M$ with boundary $\NN=\partial M$
the tangent bundle $T_*M|_{\NN}$ is trivialized
by a 4-framing $\zeta$ (usually, it is a trivialization of $T_*\NN$ stabilized by the normal vector bundle), then
the {\it  relative Pontryagin class} $p_1(M,\zeta)\in H^4(M,\NN)$ is, by definition, the obstruction class for extending $\zeta$
to $n=4$ sections of $T_*M$ without dropping the dimension below $n-1=3$. The corresponding integer {\it Pontryagin number} will be denoted
$p_1[M,\NN]$ (where $\zeta$ is supposed to be clear from the context).

If a 4-manifold $M$ is oriented
(and so is its boundary $\NN$) and $\bar M$ is obtained from it by reversing orientation (and so is $\bar\NN=\partial\bar M$),
and the framing $\bar\zeta$ is obtained by reversing one vector-component in a framing $\zeta$ on $T_M|_\NN=\nu_\NN\oplus T_*\NN$
(where, following common boundary orientation convention, $\nu_\NN$
is the normal bundle to the boundary with the outward orientation) then
\[p_1(M,\zeta)=p_1(\bar M,\bar\zeta).\]

For a closed manifold $Z$ endowed with a tangent framing $\zeta_Z$,
and any tangent framing of $M$ along $\partial M=\NN$, we have (also straight from the definition)
\[p_1(M\times Z,\zeta\times\zeta_Z)=p_1^*(M,\zeta),\]
where $p_1^*$ refers to the pull-back induced by the projection $M\times Z\to M$.
Applying it for a 2-torus $T$ and $(M,\NN)=(D^2,S^1)$ we may conclude
 that  (since any Lie-group invariant framing $\eta$ on $S^1\times T$
is isotopic to the product-framing)
\begin{equation}\label{eta-vanishing}
p_1(D^2\times T,\eta)=0.
\end{equation}

Due to another straightforward calculation,

for the normal framings $\hat\xi_n, n\ge 4,$ of $S^3$ we have
\begin{equation}\label{signed-formula}
 p_1[D^4,S^3]=-2.
\end{equation}

If we consider another compact smooth oriented 4-manifold $M'$
with $\partial M'=\bar\NN$, then the {\it Pontryagin additivity}
is the usual obstruction splitting
$p_1(M,\zeta)\oplus p_1(M',\bar\zeta)\mapsto p_1(M\cup M')$
induced by the homomorphism
$H^4(M,\NN)\oplus H^4(M',\bar\NN)\to H^4(M\cup M')$.
This gives
\be\label{P-additivity}
p_1[M,\NN]+p_1[M',\bar\NN]=p_1[M\cup M'].
\ee

\begin{lemma}\label{point-obstruction} 
Assume that a framing $\zeta$ of the tangent bundle of $D^4$ is defined
everywhere except one interior point.
Then the normal n-dimensional framing, $n\ge5$, on $S^3$ matching $\zeta|_{S^3}$
represents a multiple $kg_n\in \Pi_n^3$, $k\in\Z$, of $g_n$
with $k=\frac12 p_1[D^4,S^3]$ relative to $\zeta|_{S^3}$.
\end{lemma}

\begin{proof} It follows from Lemma \ref{spherical-classes}
and formula (\ref{signed-formula}).
\end{proof}

\subsection{Characteristic tori}\label{char_tori}
Let $M$ be a smooth closed oriented 4-manifold and
$T\subset M$ a torus with $T\circ T=0$. By a {\it torus surgery along $T$}, or shorter {\it T-surgery},
we mean a surgery consisting in cutting out from $M$ a tubular neighborhood $N(T)\cong T\times D^2$ of $T$ and gluing it back
to the compact complement,
$\M= M\sm \operatorname{Int} N(T)$, along some
orientation preserving
diffeomorphism $\phi:T\times S^1\to T\times S^1$ to obtain $M_\phi=\M\cup_\phi (T\times D^2)$.
The circles
$\phi(\pt\times S^1)\subset \partial N(T)=\partial \M$ (which are bounding the discs $\phi(\pt\times D^2)$ in $M_\phi$)
are called {\it meridians} of the $T$-surgery.

Note that the class
$w_2(\M,\eta)$
is the coboundary $\delta\a$ for some class $\a\in H^1(\partial N(T);\Z/2)$, since the absolute class
$w_2(\M)$ vanishes.

\begin{lemma}\label{char-torus}
Assume that $w_2(M)\ne0$ and a torus $T\subset M$, $T\circ T=0$, is characteristic.
Then a $T$-surgery produces $M_\phi$ with $w_2=0$, or otherwise $T$ is characteristic in $M_\phi$. More precisely,
\begin{enumerate}\item
$w_2(M_\phi)=0$ if and only if there exists
$\a\in H^1(\partial N(T);\Z/2)$ with $w_2(\M,\eta)=\delta\a$ which
takes value $1$ on the meridians of the $T$-surgery;
\item
torus $T$ remains characteristic in $M_\phi$ if and only if
there exists
$\a\in H^1(\partial N(T);\Z/2)$ with $w_2(\M,\eta)=\delta\a$ which takes value $0$ on the meridians of the $T$-surgery.
\end{enumerate}
\end{lemma}

\begin{proof}
The second homomorphism in the Mayer-Vietoris sequence (with $\Z/2$-coefficients)
\[
H^1(\partial \M)\to
 H^2(\M,\partial \M) \oplus H^2(N(T),\partial N(T))\to H^2(M_\phi)
\]
decomposes class $w_2(M_\phi)$ as the image of
$w_2(\M,\eta)\oplus w_2(N(T),\eta)$.
 So, it vanishes if and only if both relative classes are congruent modulo the coboundary homomorphisms
$ \delta_{\M} :  H^1(\partial \M)\to  H^2(\M,\partial \M)$ and
$ \delta_{N(T)} \circ \phi^* : H^1(\partial \M)\to  H^2(N(T),\partial N(T))$.
On the other hand,
(\ref{torus-w2}) implies that
$w_2(N(T),\eta)=(\delta_{N(T)}\circ \phi^*)\a$ if and only if
 the class $\phi^*\a$ takes value 1 on $\pt\times S^1$, and thus if and only if
 $\alpha$ takes value 1 on the meridians $\phi(\pt\times S^1)$ of the $T$-surgery.
 This proves item (1).

To prove item (2), we note that $T$ remains characteristic in $M_\phi$ if and only if
the image of
$w_2(\M,\eta)\oplus 0$ is zero
and apply once more (\ref{torus-w2}).
\end{proof}

We say that a characteristic torus $T\subset M$ is {\it odd} if
$w_2(M)\ne 0$, $T\circ T=0$, and $T$ remains
characteristic
under $T$-surgery for any orientation preserving diffeomorphism $\phi : T\times S^1\to T\times S^1$.

\begin{proposition}\label{essential-torus}
Let $T\subset M$ be a characteristic torus with $T\circ T=0$ and $N(T)\cong T\times D^2$ is its tubular neighborhood with the compact complement
$\M=M\sm\operatorname{Int} N(T)$.
If $w_2(M)\ne 0$, then
the following conditions are equivalent:
\begin{enumerate}\item
Torus $T$ is odd.
\item
$w_2(\M,\eta)=0$.
\item
There exists an extension of $\eta$ from the boundary to all but a finite number of points of
$\M$.
\end{enumerate}
\end{proposition}

\begin{proof}
Equivalence of (1) with (2) follows from Lemma \ref{char-torus}.
Property (3) is equivalent to (2), since
$w_2(\M,\eta)$ is the obstruction to extend
$\eta$ to the 2-skeleton of
$\M$ and its further extension to the 3-skeleton is unobstructed due to $\pi_2(\SO(4))=0$.
\end{proof}

\begin{lemma}\label{torus-framing}
Suppose that $T\subset M$, $T\circ T=0$, is an odd characteristic torus $($in particular, $w_2(M)\ne0${\rm)},
and the inclusion homomorphism $H_1(T;\Z/2)\to H_1(M;\Z/2)$ vanishes.
Consider a framing $\zeta$ in
the 3-torus $T^3=\partial N(T)$.
Then  $\zeta$ is homotopic to $\eta$
in the complement of a point of $\partial N(T)$ $($or equivalently, on its 2-skeleton$)$ if and only if
$\zeta$ augmented with an outward-normal framing
extends to all but finitely many points of $\M=M\sm \operatorname{Int} N(T)$.
\end{lemma}

\begin{proof}
Like in the above lemma,
a criterion for existence of such an extension is
$w_2(\M,\zeta)=0$.
The first obstruction to a homotopy between $\zeta$ and $\eta$
 is a class $o_1(\zeta,\eta)\in H^1(T^3;\pi_1(\SO(4)))=H^1(T^3;\Z/2)$
and our claim follows from $w_2(
\M,\zeta)-w_2(\M,\eta)=\delta o_1(\zeta,\eta)$. Namely, our
assumptions imply that $H^1(
\M;\Z/2)\to H^1(\partial\M)$ vanishes (since every 1-cycle in $\partial
\M$ bounds in $\M$) and thus, $\delta$ is monomorphic.
\end{proof}

\begin{lemma}\label{generating}
Let $\zeta$ be a framing along a 3-torus $T^3$ in its stabilized
tangent bundle $T_*(T^3\times\R)$,
and $\hat\zeta$
be a normal $n$-dimensional framing matching
 $\zeta$, $n\ge5$. Then $\hat\zeta$
represents an element $kh_n+lg_n\in\Pi^3_n$, for some $k\in \{0,1\},l\in\Z$.
\end{lemma}

\begin{proof}
Since $g_n=J(b_n)$ where $b_n$ is a generator of $\pi_3(SO(n))$, we may assume (replacing $\zeta$ by $\zeta+\xi$, if needed)
that $\zeta$ is a stabilization of a tangent framing of $T^3$.

If $\zeta$ becomes homotopic to $\eta$ after the 3-torus $T^3 $ is punctured at a point, then $\hat\zeta$ becomes also homotopic
to $\hat\eta$ over the punctured torus, and the latter homotopy produces a normal $n$-dimensional framing  on the cylinder $T^3\times[0,1]\subset \R^{n+3}\times\R$ connecting $\hat\eta$ on $T^3\times0$ and $\hat\zeta$ on $T^3\times1$
and defined at all but one interior point of the cylinder.
This shows that the element of $\Pi^3_n$ determined by $\hat\zeta$ can be expressed as $h_n+J(x)$ for some $x\in\pi_3(\SO(n))$,
and it is left to use Lemma \ref{spherical-classes}.

If $\zeta$ is not homotopic to $\eta$ over the punctured torus $T^3$,
then  it is not homotopic to $\eta$ over its 1-skeleton (since $\pi_2(SO(n))=0$), so,
there exists a non-zero element in $H_1(T^3;\Z/2)$ such that, along any loop in $T^3$ representing this element,
$\zeta$ differs from $\eta$
by a full twist (via a homotopically non-trivial loop in $\SO(3)$). We peak such a loop
representing a primitive element in $H_1(T^3;\Z)$
and consider a filling $M=T^2\times D^2$ of $T^3$ for which this loop becomes a meridian.
Then $w_2(M,\eta)\ne0$ (see (\ref{torus-w2})) implies that
$w_2(M,\zeta)=0$, which means that $\zeta$ can be extended from $T^3$ to a punctured $M$.
Hence, like in the previous case, the element of  $\Pi^3_n$ determined by
$\hat\zeta$ belongs to the image of $J$-homomorphism, and it is left again to use Lemma \ref{spherical-classes}.
\end{proof}

\subsection{Membranes to recognize odd tori}\label{membranes-section}
Let $M$ be a smooth 4-manifold.
Assume that $\S\subset M$ is the image of a compact surface smoothly immersed in
$M$ and embedded near the boundary $\partial \S$.
Suppose that $\upsilon$ is a vector field, which is
non-zero and normal to $\S$ along $\partial\S$.
Then one can define the index $\ind(\S,\upsilon)\in \Z/2$
of $\upsilon$ as the
self-intersection number of $\S$ as it is shifted by a vector field which coincides with $\upsilon$ on $\partial\S$.
Clearly, it does not change under
isotopy of $\S$ together with a continuous deformation of the non-zero normal field $\upsilon$.

Furthermore, we may
cut $\S$ into several pieces, $\S_1,\dots,\S_k$,
along some embedded curve $L\subset \S$, $\partial L\subset\partial\S$ which
does not pass through the self-intersection points of $\S$.
Let $\til\upsilon$ denote an extension of
$\upsilon$ from $\partial \S$ to $L\cup \partial \S$, still everywhere non-zero and normal to $\S$.
Then
\be\label{index-additivity}
\ind(\S,\upsilon)=\ind(\S_1,\til\upsilon)+\dots+\ind(\S_k,\til\upsilon).
\ee

For any closed orientable surface $F\subset M$,
{\it a membrane on $F$} is the image of
an immersed compact surface $\Sigma\subset M$
which approaches $F$ along its non-empty
boundary, $C=\partial\Sigma=\S\cap F$,
near which $\S$ is supposed to be embedded and nowhere tangent to $F$.
A membrane will be called {\it simple} if $\partial\S$ is connected.
The {\it index}\, $\Ind(\S)\in \Z/2$
of a membrane $\S$ is the self-intersection index $\ind(\S,\nu_C)$
of $\S$ relative to a non-zero vector field $\nu_C$ on $C$ which is normal to $C$ in $F$.
This index does not depend on the choice of $\nu_C$.

\begin{lemma}\label{membranes}
If $F$ is an orientable characteristic  surface
in a compact orientable 4-manifold $M$, and $\S$ is a membrane on $F$
then $\Ind(\S)$ depends only on the homology class of $C=\partial\S$ in $H_1(F;\Z/2)$.
In particular, $\Ind(\S)=0$ if $[C]=0$.
\end{lemma}

\begin{proof}
Consider two membranes $\S_1$ and $\S_2$, whose boundaries,
$C_i=\S_i\cap F$, are $\Z/2$-homologous in $F$.
We assume first that $C_1$ and $C_2$ are disjoint. The homology between $C_1$ and $C_2$ is realized
by a surface $R\subset F$ bounded by $C_1\cup C_2$ and the union $\hat R=R\cup\S_1\cup\S_2$ is a closed surface piecewise-smoothly immersed in $M$.
Consider along $C_i$, $i=1,2$, a non-zero vector field $\nu_i$ orthogonal
both to $F$ and to $\S_i$.
We extend $\nu_i$ to generic normal fields on $\S_i$ (with a finite set of zeros), and also
extend $\nu_1$, $\nu_2$ to a generic normal vector field on $R$. As a result we obtain some vector field  $\nu$ on $\hat R$.
Let $\hat R'$ (respectively, $R', \S_1', \S_2'$) be a small shift of $\hat R$ (respectively, $R,\S_1,\S_2$) in the direction of $\nu$. Then,
\[
\hat R\circ\hat R'=R\circ R'+\S_1\circ\S_1'+\S_2\circ\S_2'=\hat R\circ F+\Ind(\S_1)+\Ind(\S_2)
\]
which together with the Wu formula $[\hat R]^2=\hat R\circ F$
gives  $\Ind(\S_1)=\Ind(\S_2)\in \Z/2$.
In particular, this gives $\Ind(\S)=0$ for any membrane $\S$ such that $[\partial \S]=0\in H_1(F;\Z/2)$.

In the case $C_1\cap C_2\ne\varnothing$, by an isotopy of membranes we make
$\S_1$ and $\S_2$ intersect transversally and consider non-zero normal to $C_i$, $i=1,2$,
tangent to $F$ vector fields $\nu_{C_i}$.
Next,
we transform the union $\S_1\cup\S_2$ into one membrane, $\S$, by inserting at each point $p$ of $C_1\cap C_2$
a standard
smoothing of nodal boundary singularity.\footnote{
That is the smoothing which is modeled on "a half" of the standard real smoothing
of a plane nodal singularity $xy=0$ in $\C^2$: it
replaces the union of halves of coordinate lines, $\{ Im \, x\,\ge0, y=0\}\cup \{Im\, y\le0, x=0\}$,
by a half of a real conic, $\{xy=\varepsilon, Im\, x\ge 0\} = \{xy=\varepsilon, Im\, y\le 0\}$, where $0\!<\!\varepsilon<\!\!<\!1$.}
This smoothing can be performed so that the vector fields $\nu_{C_1}$ and $\nu_{C_2}$ taken
outside a small neighborhood of $C_1\cap C_2$
can be extended to a non-zero, normal to $C=\partial \S$ in $F$, vector field $\nu_C$.
Then, applying the additivity property
(\ref{index-additivity}) we conclude that
the difference $\Ind(S)-\Ind(\S_1)-\Ind(\S_2)$ is a sum of equal contributions, $\delta$, given
by the smoothing at each of the points of $C_1\cap C_2$, or equivalently,
\be\label{quadratic-property}
\Ind(\S)=\Ind(\S_1)+\Ind(\S_2)+ \delta \card(C_1\cap C_2).
\ee

On the other hand, $[\partial\S]=[\partial \S_1]+[\partial \S_2]=0\in H_1(F;\Z/2)$ and, hence, as it is proved above, $\Ind(\S)=0$.
Now, it remains to notice that $\card(C_1\cap C_2)$ is even, since $[C_1]=[C_2]\in  H_1(F;\Z/2)$ and $F$ is orientable. \end{proof}

It will be convenient to relate the index  $\Ind(\S)$ to another one.
Namely, let us consider a compact oriented 4-manifold $M'$
whose boundary $\NN'=\partial M'$
is endowed with a framing $\zeta$ of $T_*(\NN')$,
and a compact {\it immersed pair}
$(\S',C')\subset (M',\NN')$, that is the image $\S'\subset M'$ of an immersed surface, which is an embedding near its boundary
$C'=\partial\S'=\S'\cap N'$.
Then, one can define
the index $\ind(\S',\zeta)$ to be the self-intersection number $\S'\circ\S'\mod2$ as $C'$ is shifted in $\NN'$ by the normal vector field along $C'$, which
is compatible with $\zeta$ (that is together with the tangent to $C'$ vector field extends to a framing of $T_*(\NN')$ along $C'$ which is related to $\zeta$ by a null-homotopic
loop in $\SO(3)$).

\begin{lemma}\label{adjunction}
Assume that $\partial M'=T^3$ and
$(\S',C')\subset (M',T^3)$ is an immersed pair with
$C'$ connected.
Then
\[
w_2(M',\eta)[\S',C']=\ind(\S',\eta)+1
\]
\end{lemma}

\begin{proof} By adjunction,
$w_2(M',\eta)[\S', C']= \ind(\S',\eta) + (w_2(\S',C')+ w^2_1(\S',C'))[\S',C']=\ind(\S',\eta)+1$,
where $(w_2(\S',C')+ w^2_1(\S',C'))[\S',C']=1$
(for closed surfaces, $w_2+w_1^2=0$, and in the case of one boundary component, $\partial S'=C'$,
$0$ is replaced by $1$).
\end{proof}

If $w_2(M')=0$, then $\ind(\S',\zeta)\in\Z/2$ depends only on the class $[C']\in H_1(\NN';\Z/2)$ and the homotopy class of $\zeta$.
So, in such a case we may skip $\S'$
from the notation letting $\ind(C',\zeta)=\ind(\S',\zeta)$ and viewing it as a $\Z/2$-valued function
(depending only on the homotopy class of $\zeta$)
defined on the kernel of
$H_1(\NN';\Z/2)\to H_1(M';\Z/2)$.

\begin{cor}\label{constant-criterion}
Let $M$ be a closed oriented 4-manifold with $w_2(M)\ne0$, $T\subset M$ be a characteristic torus with
$T\circ T=0$ that has a trivial inclusion homomorphism $H_1(T;\Z/2)\to H_1(M;\Z/2)$. Let $N(T)$ be a small tubular neighborhood of $T$ in $M$,
and $\M$ the compact manifold obtained by removing $N(T)$ from $M$.
Then $T$ is odd if and only if
$\ind(C',\eta)$
as a function
$H_1(\partial \M;\Z/2)= H_1(\partial N(T);\Z/2)  \to\Z/2$
takes value $1$ on each non-zero element.
\end{cor}

\begin{proof}
It follows immediately
from Lemmas \ref{essential-torus} and \ref{adjunction} (the latter is applied to $M'=\M$).
\end{proof}

Now, assume that $M$, $T$, $N(T)$, and $\M$ are as in Corollary \ref{constant-criterion} above, and consider a {\sl simple membrane}
$\S\subset M$ on $C=\partial\S\subset T$.

The surface $\S'=\S\cap \M$
with $C'=\partial\S'$ gives an
immersed pair $(\S',C')\subset(\M,\partial \M)$
and we may consider $\ind(\S',\eta)$ with respect to a
Lie-group invariant  tangent framing $\eta$
on $T^3=\partial \M$, which is well defined up to homotopy.

\begin{lemma}\label{two-indices}
If the class $[C]\in H_1(T;\Z/2)$ does not vanish, then $\Ind(\S)=\ind(\S',\eta)$.
If $[C]=0$, but $[C']\ne0$,
then $\Ind(\S)=\ind(\S',\eta)+1$.
\end{lemma}

\begin{proof} It follows from
an observation
that, for any smooth simple closed curve in $T$, the framing formed along the curve by a tangent and a normal vector is compatible with the Lie-group invariant tangent framing on $T$ if and only if the homology class of the curve in $H_1(T;\Z/2)$ is non-zero.
\end{proof}

\begin{proposition}\label{criterion} Let $M$ be a smooth closed oriented 4-manifold with a torus
$T\subset M$ realizing the class $w_2(M)\ne 0$, so that $T\circ T=0$ and
the inclusion homomorphism $H_1(T;\Z/2)\to H_1(M;\Z/2)$ is trivial.
Then $T$ is an odd torus if
for any non-trivial class in $H_1(T;\Z/2)$
there exists a simple membrane
of odd index whose boundary realizes this class.
\end{proposition}

\begin{proof}
It follows immediately from Lemma \ref{two-indices} and Corollary
\ref{constant-criterion}.
\end{proof}

\subsection{Divisibility of {\boldmath$\sigma$} and the order of \boldmath$\Pi_n^3$}
In what follows $M$ is a compact oriented 4-manifold with boundary $\NN=\partial M$ and $\zeta$ is
a framing of $T_*M|_\NN$. We can extend $\zeta$ to a framing
$\zeta_M$ of $T_*M$ defined over the complement of a closed characteristic surface $F\subset M$
and a few points.

\begin{lemma}\label{kill-sphere}
Assume that $w_2(M,\zeta)\in H^2(M,\NN;\Z/2)$ is non-zero and realized by a sphere
$F\subset M\sm \NN$ with $F\circ F=0$. Then a Morse
surgery of index $3$
replacing a tubular neighborhood of $N(F)\cong S^2\times D^2$ with $D^3\times S^1$ results in a manifold $M'$
with $w_2(M',\zeta)=0$.
\end{lemma}

\begin{proof}
Put $\M=M\sm \operatorname{Int}N(F)$.
The second homomorphism in the Mayer-Vietoris sequence (with $\Z/2$-coefficients)\[
H^1(S^2\times S^1)\to
 H^2(\M,L\cup (S^2\times S^1)) \oplus H^2(D^3\times S^1,S^2\times S^1)\to H^2(M',L)
\]
decomposes class $w_2(M',\zeta)$ as the image of
$w_2(\M,\zeta_M\vert_{\partial \M})\oplus w_2(D^3\times S^1,\zeta_M\vert_{\partial (S^2\times S^1)=\partial N(T)})$. The first summand is zero, since the framing $\zeta_M$
is well defined at all but finite number of points in $\M$. The second one is zero, since
$ H^2(D^3\times S^1,S^2\times S^1)=H_2(D^3\times S^1)=0$.
\end{proof}

\begin{lemma}\label{cutting}
Let $F\subset M$ be a closed orientable surface which is characteristic for $(M,\zeta)$ and
$\Sigma$ a simple membrane
of even index on $F$.
Consider a manifold $M'$ that is obtained by
a Morse surgery of index 2
on $M$ which cuts out a neighborhood $N(C)\cong S^1\times D^3$ of $C$ and fills the complement with $D^2\times S^2$.
If the diffeomorphism  $N(C)\cong S^1\times D^3$
is chosen to respect
the framing of $C$ defined by two vector fields,
 normal to $C$ in $\S$ and normal to $C$ in $F$, then
the band $F\cap N(C)\cong S^1\times[0,1]$ can be replaced by a pair of disjoint discs inside $D^2\times S^2$ to form a surface $F'\subset M'$
which is characteristic for $(M',\zeta)$.
\end{lemma}

\begin{proof}
By an isotopy in $M$ we identify $(N(C), C)$ with $(N(C'), C')$ where $C'\subset \Sigma$ is obtained by pushing $C$ inside $\Sigma$.
Thus,  $M'$ can be seen as a Morse surgery along $C'$. We put $\M=M\sm \operatorname{Int}N(C')$ and notice that the former surface $F'$ is $\Z/2$ homologous to $F$, which in this new model of $M'$
is no more affected by the surgery.
The Mayer-Vietoris sequence (with $\Z/2$-coefficients) written for this model,
\[
H^1(S^1\times S^2)\to
 H^2(\M,L\cup (S^1\times S^2)) \oplus H^2(D^2\times S^2,S^1\times S^2)\to H^2(M',L),
\]
decomposes class $w_2(M',\zeta)$ as the image of
$w_2(\M,\zeta_M\vert_{\partial \M})\oplus w_2(D^2\times S^2,\zeta_M\vert_{S^1\times S^2})$. Due to the initial choice of the framing $\zeta_M$, the first summand is Poincar\'e dual to $F$.
The second summand is either $0$ or Poincar\'e dual to the core-sphere $S= 0\times S^2$. To exclude the second option, we consider the sub-surface $\Sigma'\subset \Sigma$ bounding $C'$,
complete it up to a closed surface $\hat \Sigma$ by a disc $D^2\times \pt$, and observe that $\hat\Sigma\circ \hat\Sigma= 0\mod 2$ while $\hat\Sigma\circ (F\cup S)=1$,
which contradicts to (\ref{Wu-relation})
\end{proof}

\begin{lemma}\label{cut-into-tori} If $H_1(M)=0$, then
using some surgery operations on $M\-\NN$,
one can obtain another compact oriented 4-manifold $M'$, $\partial M'=\NN$,
with a characteristic surface $F'\subset M'$ formed by several disjoint
tori (possibly, $F'=\varnothing$).
\end{lemma}

\begin{proof}
The condition  $H_1(M)=0$ implies that both $H^2(M)$ and $H^2(M,\NN)$ are torsion free, and so, class $w_2(M,\zeta)\in
H^2(M,\NN; \Z/2)$ is integral
and can be represented by an orientable closed connected surface $F$.
If $F$ has genus $0$, then by taking connected sum with
$\barCp2$ or
$\Cp2$ at the points of $F$, we achieve $F\circ F=0$ and still
keep $F$ being a characteristic surface in the resulting 4-manifold. Then we apply Lemma \ref{kill-sphere} and obtain $M'$ with $F'=\varnothing$.

If $F$ has genus $g>1$, it contains a disjoint set of simple closed
curves $C_1,\dots,C_{g-1}$ which cut $F$ into $g$ components, which are tori with holes.
Each curve $C_i$ is null-homologous in $F$, which implies that it bounds a membrane
 of index $0$ by Lemma \ref{membranes}.
Applying Morse surgery of Lemma \ref{cutting}, we obtain a 4-manifold $M''$ where the characteristic surface $F''$ splits into $g$ closed torus components.
By several blowups (adding $\barCp2$) and anti-blowups (adding $\Cp2$) at the points of $F''$ we obtain manifold $\hat M''$ in which the obtained surface $\hat F''$
is still characteristic, with torus components having
trivial normal bundles.
\end{proof}

Now we can prove the basic relations between
the elements $g_n, h_n\in\Pi^3_n$ ({\it cf.}, \cite{R1,R2}).

\begin{proposition}\label{generate}
$g_n$ and $h_n$ generate $\Pi^3_n$ for $n\ge 5$. The order of $g_n$ is the minimal nonzero value of $|\frac{p_1[M]}2|$ over all smooth closed oriented $4$-manifolds $M$ with
$w_2=0$.
\end{proposition}

\begin{proof}
Pontryagin construction gives a representation of any given element of  $\Pi^3_n$ by
a closed oriented smooth 3-manifold $\NN$ with a normal $n$-dimensional framing. It corresponds to some framing,
$\zeta$, in the stabilized tangent space of $\NN$.
Due to Rokhlin's Cobordism Theorem \ref{cobordism},
we can find an oriented 4-manifold $M$ bounding $\NN$. Moreover, $M$ can be assumed to be simply-connected (which is achieved by Morse surgery of index 2).
Applying Lemma \ref{cut-into-tori}, we obtain manifold $M'$ whose characteristic surface $F'$ splits into tori (or empty).

This means that $\zeta$ can be extended to $M'$ in the complement of $F'$
and several punctures. This complement gives a framed cobordism between $(L,\zeta)$ and the union of
several framed 3-tori with several framed 3-spheres.

By Lemma \ref{generating} each framed 3-torus represents an element $kh_n+lg_n\in
\Pi_n^3$, $n\ge5$, for some $k\in\{0,1\}, l\in\Z$,
while the framings on the 3-spheres represent by Lemma \ref{point-obstruction}
multiples of $g_n$. This proves the first claim of the proposition.

As for the order of $g_n$, note that, according to Pontryagin construction, $l g_n=0$ with $l> 0$ if and only if $l$ copies of $S^3\subset \R^{n+3}$ with a normal framing $\hat\xi_n$
bound a normally framed 4-manifold, $(\M, \xi)$ in $\R^{n+4}$.
 After filling the boundary spheres of $\M$ with 4-discs, we obtain a closed 4-manifold $M\subset \R^{n+5}$ with $w_2=0$.
Then, due to
Pontryagin's additivity (\ref{P-additivity})
and formula (\ref{signed-formula}), we get
$$p_1[ M]=-p_1(\nu_M)[M]=-p_1(\M,\hat\xi_n)[\M,\partial\M] -l p_1(D^4,\hat\xi_n)[ D^4, S^3]= 2l$$
(here $\nu_M$ stands for the normal bundle of $M$).
Reciprocally, given a closed 4-manifold $ M$ with $w_2=0$, we consider its tangent framing $\zeta$ defined at all but finite number of points in $M$.
According to Lemma \ref{point-obstruction} and the additivity relation (\ref{P-additivity}),
such a framing represents a relation $lg_n= 0$ with $l=\sum l_i\in \Z$ where the sum is taken over
the singular points of the framing $\zeta$ and $l_i = \frac12 p_1[D^4_i,S^3_i]$ for an $i$-th singular point with $p_1$ taken relative to $\zeta|_{S^3}$, so that $l=\sum l_i=\sum \frac12
p_1[D^4_i,S^3_i]= \frac12 p_1[M]$. (Note that in the last argument the sign of $l$ comes out the same as the sign of $p_1[M]$.)
\end{proof}

\begin{proposition}\label{minimal} For $n\ge 5$,
we have
$12g_n=h_n$.
\end{proposition}

\begin{proof}
Let $M=\Cp2\#9\barCp2$ be obtained by blowing up 9 points in a general position in $\Cp2$, and let $T\subset M$ be a torus given
given by a nonsingular plane cubic curve traced through the chosen 9 points and then lifted to $M$. Such torus $T$ has zero self-intersection and realizes $w_2(M)$.
Furthermore, $H_1(T;\Z/2)$ is generated by vanishing circles, which bound
disc membranes of index $-1$, so that Proposition \ref{criterion} applies and shows that
$T$ ia an odd characteristic torus.

Next, we consider a tubular neighborhood $N(T)=T\times D^2$ of $T$ in $M$ and a framing of $T_*M$ along $T^3=\partial N(T)$ given by $\eta$ augmented by a normal vector field.
By Lemma  \ref{torus-framing}, this framing extends to a tangent framing of $\M=M\sm \operatorname{Int}N(T)$ except a finite set of points  $x_1,\dots, x_k\in \M$.
We denote this extended framing by $\zeta$. Its restriction to the complement $M'= \M\sm \bigcup_i D_i^4$ in $\M$ of the union of small balls $D_i^4$ around the points
 $x_i$, $i=1,\dots,k$, provides a framed cobordism  leading to a relation $h_n=l g_n$ for some $l\in\Z$.

Finally, we find $\l$ combining
the relation $p_1(M)=3\sigma(M)=-24$
with
the additivity relation (\ref{P-additivity}) that gives
$$
p_1[M]=\sum_i p_1[D^4_i,S^3_i]+p_1[T\times D^2, T\times S^1]+p_1[M',\partial M'].
$$
The last summand vanishes, since $\zeta$ is defined at all the points of $M'$.
The second summand
$p_1[T\times D^2, T\times S^1]$ also vanishes (see (\ref{eta-vanishing})),
while the first summand is formed by
$p_1[D^4_i,S^3_i]=2l_i$ where according to Lemma \ref{point-obstruction} $\zeta\vert_{S^3_i}$ matches $l_i\hat\xi_n$. Hence, $h_n=l g_n$ with $l=\sum l_i=-12$. The sign of $l$ can be switched to opposite due
to $2h_n=0$ (see  Lemma \ref{tori-framings}(2)).
\end{proof}

\begin{cor}\label{pontryagin48}
If $h_n\ne 0$ for $n\ge 5$, then the minimal nonzero absolute value of $\frac{p_1(M)}2$ over all smooth closed oriented $4$-manifolds $M$ with zero second Stiefel-Whitney class is equal to $24$ and, as so, $p_1(M)$ is divisible by $48$ and $\sigma(M)$ is divisible by $16$ for any such manifold $M$.
\end{cor}

\begin{proof} Due to Proposition \ref{minimal}, if $h_n\ne 0$ the order of $g_n$ is $24$, which
being combined with Proposition \ref{generate}
implies that  the minimal nonzero absolute value of $\frac{p_1(M)}2$ is at least $24$. On the other hand,
$\frac{p_1[M]}2=24$ when $M$ is a K3-surface. Since the set of values of $p_1[M]$ forms a group, it implies the divisibility of $p_1[M]$ by $48$ and, then, due to Theorem \ref{formula} the divisibility of
$\sigma(M)$ by $16$.
\end{proof}

\subsection{The final step: $\mathbf{h_n\ne0}$}\label{final-step}
According to Corollary \ref{pontryagin48}, to complete the proof Theorem \ref{divisibility} it remains to show that $h_n\ne 0$ for $n\ge 5$.

\begin{lemma}\label{condition}
If $h_n=0$ for $n\ge 5$, then there exists a smooth closed oriented 4-manifold $M$ with the following properties:
\begin{itemize}
\item[a)] $M$  contains an odd characteristic
torus $T\subset M$, $T\circ T=0$.
\item[b)] $p_1[M] =0$.
\end{itemize}
\end{lemma}

\begin{proof}
We use the same arguments as in the proof of Proposition \ref{minimal}: take the manifold from this proof
and replace the balls around the points
which obstruct the framing by a manifold with boundary to which the framing extends, which is possible if $h_n=\sum\pm g_n$ is zero.
\end{proof}

\begin{proposition}\label{inconsitancy} Conditions a) and b) are inconsistent.
\end{proposition}

\begin{proof}
The proof is by contradiction. We assume that a smooth closed oriented 4-manifold $M$ satisfies both properties, a) and b).
Note that we can also assume $M$ to be simply-connected. Indeed, we can start from considering a tangent framing in $M\- T$ as in Proposition \ref{essential-torus}(3)
and make a sequence of Morse modifications of index $2$ along circles $S^1_i$ in the complement of the singular set of the framing with a precaution to choose each time
that of two possible modifications (which differ by a choice of the system of parallels on the boundary of $S_i^1\times D^3\subset M$) to which the framing extends.
Thus, we kill the fundamental group preserving the properties a) and b).

By Theorem \ref{cobordism}(2), $M=\partial W$ for some oriented 5-manifold $W$. The class $w_2(W)$ is dual to a compact 3-manifold
$Q\subset W$ and, since $T\circ T=0$, we may assume that $\partial Q=T$ (see \cite[Commentaires sur le quatri\`eme article]{Book}
 for a detailed justification).

 We choose a tubular neighborhood $U\subset W$ of $Q$, with its $D^2$-fiber projection $U\to Q$,
 and note that the boundary $X=\partial U$ is a closed 4-manifold glued from two pieces: the total space $\X$ of
 the associated $S^1$-fiber bundle $q:\X\to Q$ and a tubular neighborhood  $V=X\cap M\cong T\times D^2$ of $T$ in $M$
(with a corner along the common boundary $T\times S^1$).

The Poincar\'e-Lefschetz duality in $(Q,T)$ implies
that some simple closed
curve $C\subset T$ realizing a non-zero element in $H_1(T;\Z/2)$
is null-homologous in $Q$ and, thus, bounds a surface $F\subset Q$.
With respect to the Lie-group invariant vector field $\eta^2$ on $T$
(see Section \ref{framed-tori}) we have the following.

\begin{lemma}\label{index-of-F}
Class $(w_2+w_1^2)(Q,\eta^2)$ evaluated on $[F,C]\in H_2(Q,T;\Z/2)$ gives $1$.
\end{lemma}

\begin{proof}
Let us fill in the boundary $T$ of $Q$
with
a solid torus $S^1\times D^2$ whose meridian is $C=\pt\times S^1$.
In the closed manifold $\hat Q$ that we obtain, the sum
$w_2+w_1^2$
vanishes. Thus, by additivity  (\ref{additivity-w})
\[(w_2+w_1^2)[F,C]=w_2[\pt\times D^2,C]+w_1^2[\pt\times D^2,C]=1+0\mod2.\hskip20mm\qed\]
\vskip-5mm\phantom\qedhere
\end{proof}

Next, we take a push-off $F'\subset X$ of $F$ by a vector field normal to $Q$ in $W$, so that $C'=\partial F'\subset \partial V=T\times S^1$
(such a vector field exists, since the surface $F'$ is not closed).

\begin{lemma}\label{index-of-F-prime}
Class $w_2(\X,\eta)$ evaluated on $[F',C']\in H_2(\X,\partial\X;\Z/2)$ gives $1$.
\end{lemma}

\begin{proof}
The tangent bundle $T_*(\X)$ splits into the pull-back of $T_*Q$ and the line bundle tangent to the $S^1$-fibers of the projection $\X\to Q$.
The latter line bundle is isomorphic to the pull-back of the determinant bundle $\det_Q$ of $T_*Q$, since $W$ is orientable.
This implies that $w_2(\X)$ is the pull-back of $w_2(T_*Q\oplus\det_Q)=w_2(Q)+w_1^2(Q)$.
The same decomposition holds for the relative classes implying that
$w_2(\X,\eta)$ evaluated on $[F',C']$ equals to
$(w_2+w_1^2)(Q,\eta^2)$ evaluated on
$[F,C]$, and the latter gives $1$ due to Lemma \ref{index-of-F}.
\end{proof}

Let us put $\accentset{\circ}W=W\sm U$ and note that
vanishing of the absolute class
$w_2(\accentset{\circ}W)$ implies that the relative one,
$w_2(\accentset{\circ}W,\eta)$,
with respect to the framing $\eta$ on the 3-torus $\partial N(T)\subset \accentset{\circ}W$
(extended by outward-normal vector fields on $\X$ and $\accentset{\circ}M$)
is a coboundary
$w_2(\accentset{\circ}W,\eta)=\delta\a$ for some
$\a\in H^1(\partial N(T);\Z/2)$.
Functoriality implies that $w_2(\X,\eta)$ and $w_2(\accentset{\circ}M,\eta)$ are restrictions of $w_2(\accentset{\circ}W,\eta)$ to
$\X$ and $\M$ respectively,
 and, thus, by Lemma \ref{index-of-F-prime}
$$
\a[C']=w_2(\accentset{\circ}W,\eta)[F',C']=w_2(\X,\eta)[F',C']=1. $$
On the other hand, $H_1(M;\Z/2)=0$
and $w_2(M)\ne 0$ imply
existence of a membrane
$\S'\subset M$, $\partial\S'=C'$, and
Proposition \ref{essential-torus} yields
$\a[C']=w_2(\M,\eta)[\S',C']=0$.
This contradiction completes the proof.
\end{proof}

\subsection{Calculation of the stable homotopy groups $\pi_{n+3}(S^n)$}\label{conclusion}
The above proof of Rokhlin's Signature Divisibility Theorem contains
at the same time a calculation of  $\pi_{n+3}(S^n)$ for $n\ge 5$.
\begin{theorem}\label{pi-stable} For each $n\ge 5$, the following holds:
\begin{enumerate}
\item
The group $\pi_{n+3}(S^n)$ is isomorphic $\Z/ 24$ and
$g_n\in \Pi^3_n=\pi_{n+3}(S^n)$ is its generator.
\item The homomorphism $J: \pi_3(\SO(n)) \to \pi_{n+3}(S^n)$ is surjective.
\end{enumerate}
\end{theorem}
\begin{proof} 
According to Propositions  \ref{generate}--\ref{minimal}, $g_n$ is a generator of $\Pi^3_n=\pi_{n+3}(S^n)$. Its order is $24$, as it follows
from Proposition \ref{generate} and Corollary \ref{pontryagin48} (this Corollary can be applied, since $h_n\ne 0$ as it is shown in Section \ref{final-step}). By definition, $g_n$ belongs to the image of $J$, which implies surjectivity of the latter.
\end{proof}

%% file: add.tex
\section{Appendix}
\setlength\epigraphwidth{.45\textwidth}
\epigraph{“The seriousness of a theorem, of course, does not lie in its consequences, which are merely the evidence for its seriousness.”}
{G.H. Hardy,  {\it A Mathematician's Apology.}}

―

\subsection{Groups $\pi_{n+3}(S^n)$ with $n\le 4$.}
By the time of Rokhlin's research, the groups $\pi_{n+3}(S^n)$ was known only for $n$ equal $1$ and $2$: $\pi_{4}(S^1)=0$ and $\pi_5(S^2)=\Z/2$ (where the latter followed from Pontryagin's calculation of $\pi_{n+2}(S^n)$ and the Hurewicz exact sequence applied to the Hopf fibtation).
The groups $\pi_6(S^3)$ and $\pi_7(S^4)$ were found by Rokhlin in \cite{R4} (after correcting
a mistake made in \cite{R2})
as a by-product of his calculation of $\pi_{n+3}(S^n)=\Z/24$
for $n\ge 5$ (independently, and using another method, the same answers were obtained
by Massey-Whitehead and Serre). In terms of generators
introduced in  Subsection \ref{framed-spheres}, and using Rokhlin's notation $\pi_{n+3}^0(S^n)=J(\pi_3(SO(n)))$, his
result
looks as follows.

\begin{theorem}[\cite{R4, SerreNote}]\label{unstable}
\quad
\begin{itemize}
\item[(a)]
$\pi_6(S^3)$ is the cyclic group $\Z/12$ with $f_3$ as a generator.
\item[(b)] $\pi_7(S^4)$ is the direct sum $\Z\oplus\Z/12$ of a free cyclic group $\Z$ with generator $g_4$
and a cyclic group with generator $f_4$.
\end{itemize}
In particular, $\pi^0_{n+3}(S^n)=\pi_{n+3}(S^n)$ for $n= 3$ and $4$.
\end{theorem}

\begin{proof}
By definition, $\pi^0_6(S^3)$ is cyclic and generated by $f_3$. Besides,
as was already known by then:
\begin{enumerate}
\item
 The suspension homomorphism $\S : \pi_6(S^3)\to \pi_7(S^4)$ is injective.
 \item  $\pi_7(S^4)=\S(\pi_6(S^3))\oplus\Z$ and
$\pi^0_7(S^4)=\S(\pi^0_6(S^3))\oplus\Z$ where the second summands are generated by $g_4$.
\item
The kernel of $\S: \pi_7(S^4)\to \pi_8(S^5)$ and that of $\S: \pi^0_7(S^4)\to \pi^0_8(S^5)$ are generated by $f_4-2g_4$
(with $f_4=\Sigma(f_3)$).
\end{enumerate}
Together with Rokhlin's theorem \ref{pi-stable} this implies the result in a straightforward way.
\end{proof}

\subsection{\boldmath$\Omega^{\Spin}_3=0$} By the time of Rokhlin's research, the notion of {\it spin structure on a manifold}
was not yet properly developed.
But nevertheless it was already present in Rokhlin's notes in a disguise of {\it a homotopy class of stable trivialization of a tangent bundle over a 2-skeleton} (usually stated as a homotopy class of stable trivializations over 1-skeleton that are extendable to 2-skeleton), which can be taken as one of many equivalent definitions of spin structure. Moreover,  Rokhlin's arguments, like those in above reconstruction of Rokhlin's proof of the Signature Divisibility Theorem, can be applied to show that his other result
$\pi^0_{n+3}(S^n)=\pi_{n+3}(S^n)$ for $n\ge 3$
implies the triviality of the spin-cobordism group in dimension 3, $\Omega^{\Spin}_3=0$ (see for details \cite{Book}).

\begin{theorem}\label{spin-zero}
Every closed oriented 3-manifold $X$ with a spin structure
bounds a smooth compact oriented spin 4-manifold $W$ whose spin structure restricts to that of $X$.
\qed
\end{theorem}

\subsection{Realization of homology classes}
The Hurewicz homomorphism theorem raised an important question which remained open for more than a decade:
{\sl Can every element of the second homotopy group of a simply connected 4-dimensional smooth manifold be represented by a smoothly
embedded sphere?} Rokhlin
answered to it negatively.

\begin{theorem}\label{non-realizability}
The class $3H\in H_2(\Cp2)$ is not realizable by a smoothly embedded sphere. $($Here $H$
stands for a generator of $H_2(\Cp2)$.$)$
\end{theorem}

\begin{proof}
Assume that a sphere $S\subset\Cp2$ realizes class $3H$. Then, $M=\Cp2\#\barCp2$
satisfies the relation $p_1(M)=0$ and contains a characteristic odd torus $T\# S$,
with zero self-intersection number, obtained by internal connected sum of
a non-singular cubic curve (2-torus) $T\subset\Cp2$ and $S\subset\barCp2$, and, so,
this contradicts to Proposition \ref{inconsitancy}.
\end{proof}

A common confusion about this discovery is due to an error in the first
Rokhlin's note on this subject, \cite{R1},
where Rokhlin claimed the opposite:
he believed that there exists an embedding of a sphere realizing class $3H$ and, thus, by Lemma \ref{condition}, $h_n=0$.
Then, Propositions \ref{generate} and \ref{minimal}
leaded him in \cite{R2} to a mistaken conclusion that $\Pi_n^3=\pi_{n+3}(S^n)=\Z/12$, which was finally corrected in \cite{R4}.

It seems that this Rokhlin's counter-example remained unnoticed
almost for a decade, until Rokhlin communicated it
to Kervaire and Milnor, who then upgraded it to a
``generalized version of Rokhlin's Signature Divisibility Theorem''.

\begin{theorem}[\cite{KM}] Let $x\in H_2(M;\Z)$  be dual to the Stiefel-Whitney class $w_2(M)$. If $x$
can be represented by a differentiably embedded 2-sphere in $M$, then the self-intersection number $x^2$ must be congruent to $\sigma(M)$ modulo 16.
\end{theorem}

\begin{proof}
By changing the orientation of $M$ to make $x^2$ nonpositive, and taking connected sum of $M$ with $s=-x^2+1$ copies of $\C P^2$, we obtain a new manifold $M_1$ and a smoothly embedded in it sphere realizing a homology class $y$ with $y^2=1$ and $y\mod 2=w_2(M_1)$. Now, there remain to notice that this embedded sphere can be blown down, which gives us a new manifold
$M_2$ with $w_2(M_2)=0$ and $\sigma(M_2)=\sigma(M_1)-1=\sigma(M)-x^2$, and then to apply Rokhlin's Signature Divisibility Theorem \ref{divisibility}.
\end{proof}

\subsection{Refinement of Rokhlin's Signature Divisibility Theorem}\label{refinement}
In the '60s, getting interested in tracking
realizability of $\Z$-homology classes in 4-manifolds by smooth surfaces of specified genus,
 Rokhlin came to a remarkable extension of the signature congruence to
4-manifolds without restriction on the
genus of a surface representing the second Stiefel-Whitney class.
Rokhlin's personal archive (that unfortunately became untraceable after the death of his wife in 1993) contained notes dated 1964 with a corresponding statement
(Theorem \ref{refined} below) and its proof.
As Rokhlin told to us, he reported this result
at the Moscow ICM-1966 at one of satellite seminars.

\begin{theorem}[Rokhlin's Refined Congruence]\label{refined} Let $M$ be a closed
oriented smooth 4-manifold and $F\subset M$ an orientable smooth characteristic surface such that the inclusion homomorphism $H_1(F;\Z)\to H_1(M;\Z)$ is zero.
Then,
$$\sigma(M)=F\cdot F + 8\, {\rm Arf}   (q _F) \mod 16$$
where ${\rm Arf}   (q _F)$  is the Arf-invariant of the quadratic function $H_1(F;\Z/2)\to \Z/2$ defined by the index of membranes.
\footnote{
The function $q_F$ is well defined due to  Lemma \ref{membranes},
and its quadratic property, $q(x+y)=q(x)+q(y)+x\circ y$,  follows from (\ref{quadratic-property}), where
$\delta=1\in\Z/2$. The latter value
is obtained immediately from consideration of the
 halves of two line-generators of $\C P^1\times \C P^1$ as
membranes on the characteristic torus $\R P^1\times \R P^1$.}
\end{theorem}

Rokhlin published this result only in 1972,
after he found its application to
a challenging Gudkov's conjecture from topology of real plane algebraic curves (see \cite{R6}).
Suggested proof of the conjecture
contained
however
a flaw which was corrected in 1977 by A.~Marin,
who in a joint work with L.~Guillou
generalized Rokhlin's refined congruence to non-orientable surfaces $F\subset M$ (see \cite{GM, Book}). In this generalization the index-function
with values in $\Z/2$ is upgraded to a function $\hat q_F : H_1(F:\Z/2)\to \Z/4$ quadratic in a sense that $\hat q_F(x+y)=\hat q_F(x)+\hat q_F(y)+2\,x\circ y$,
the Arf-invariant is replaced by the Brown-invariant, and
the congruence takes form
\be\label{GM}
\sigma(M)=F\cdot F +  2\, {\rm Br}   (\hat q _F) \mod 16.
\ee

Since then, this Guillou-Marin-Rokhlin
congruence found many other
applications not only in real algebraic geometry, but also in study of the non-oriented 4-genus of classical knots and links.


Rokhlin's proof of Theorem \ref{refined}
was based exclusively on cobordism arguments and geometrical constructions that appeared earlier in his proofs in \cite{R1} - \cite{R4}
(in particular, it provided a new, more elementary proof of Theorem \ref{divisibility} free of homotopy computations).
Regrettably, it was not published. Fortunately, almost immediately the importance of geometric content hidden behind this theorem has been recognized by
A.~Casson who unraveled its proof and included it 
into his famous (unpublished) lectures on topology of 4-manifolds in the middle of 70s.
Furthermore, Casson's proof was developed by M.~Freedman and R.~Kirby in  \cite{FK}.
It was at about the same time that L.~Guillou and A.~Marin \cite{GM}, being inspired directly by Rokhlin's paper,
obtained their generalization
of Theorem \ref{refined} to non-orientable characteristic surfaces.
One more version of the proof was obtained later by Yu.~Matsumoto \cite{Book}.
Various remakes of these cobordism-based proofs can be found in \cite{KirbyLN} and \cite{Scorpan}.

\subsection{Higher-dimensional generalizations}
One of the first higher-dimensional generalizations of Rokhlin's Signature Divisibility Theorem came from Atiyah-Hirzebruch proof of integrality conjecture for Todd genus. They designed for that a topological version of the Grothendieck-Riemann-Roch
Theorem, and then applied it successfully to a study of integrality properties of the so-called $\hat A$-genus (a special system of multiplicative polynomials introduced a little before by Hirzebruch). As one of the corollaries they obtained
the following result.

\begin{theorem}[\cite{AH}] If $M$ is a smooth closed oriented $(8n+4)$-dimensional manifold with $w_2(M)=0$,
 then its $\hat A$-genus is an even integer.
\end{theorem}
Since in dimension four $\hat A(M)=-\frac{p_1[M]}{24}$, this is a generalization of Rokhlin's divisibility theorem.

With a later discovery that $\hat A$-genus is equal to the index of the Dirac operator (see \cite{AS}),
this Atiyah-Hirzebruch proof got a very short, and elementary, modulo Atiyah-Singer Index Theorem, form. Namely,
if $w_2=0$ for a smooth closed oriented 4-manifold $M$, then a $\Spin$-structure $\phi$ in $M$ defines a pair of $\SU(2)$-bundles $V^\pm$
and the Dirac operator $\slashed{D}_\phi:\Gamma(V^+)\to\Gamma(V^-)$  has index equal to $ -p_1[M]/24$ by the Index Theorem.
The quaternionic structure in $V^\pm$ implies that this index is an even integer and it remains
to use Rokhlin's Signature Formula $p_1[M]=3\s(M)$.

However,
in dimensions higher than $4$ the signature is governed by a completely different system of multiplicative polynomials, and the question of Rokhlin-like divisibility for the signature of higher-dimensional manifolds
stayed open up to the
70th when S.~Ochanine extended Rokhlin's divisibility theorem to dimensions $> 4$.
At first, in \cite{O1}, he proved that {\sl the signature of a smooth closed oriented $(8n+4)$-manifold is divisible by 16, if the manifold admits an $SU$-structure.}
But later on, he noticed that the assumption "$SU$-structure" can be replaced by the assumption "$Spin$-structure", thus giving a perfect generalization
of Rokhlin's divisibility of signature by 16 (see \cite{Ochanine}).

He discovered also what can be a higher-dimensional analog of Rokhlin's Refined Congruence \ref{refined} by transforming the latter into framework of $\Spin^\C$-manifolds
(oriented manifolds with trivial third integral Stiefel-Whitney class $W_3$).
Namely,
he showed how an orientable characteristic submanifold $F^{8n+2}$ in a $\Spin^\C$-manifold $M^{8n+4}$ can be endowed with an induced $\Spin$-structure,
and then enriched $F$ by a certain $\Spin$-cobordism invariant  $\mathcal{K}(F)$.

\begin{theorem}[\cite{Ochanine}]
Let $M$ be smooth closed oriented
$(8n+4)$-manifold
with $W_3(M)=0$
and let $F\subset M$
its oriented characteristic submanifold. Then:
\begin{itemize}\item
$\mathcal K(F)\in\Z/2$ depends only on $F$, but not on the $\Spin^\C$-structure of $M$.
\item
$\s(M)=\s(F\circ F)+8\mathcal K(F)\mod16$, where $F\circ F$ is any of $8n$-dimensional manifolds obtained from $F$ by self-intersecting in $M$.
\end{itemize}
\end{theorem}

In 1990 S.~Finashin \cite{Finashin} extended
Ochanine's theorem to the case of a non-oriented characteristic submanifold $F\subset M$.
For such a refinement of Ochanine's theorem he replaced an auxilliary $\Spin^C$-structure on $M$ by a $\Spin$-structure on $M\-F$, endowed
$F$ with an induced $\operatorname{Pin}^-$-structure, and introduced an appropriate $\operatorname{Pin}^-$-cobordism invariant $I(F)\in\Z/8$
replacing Ochanine's $\Spin$-cobordism invariant $K(F)$.

\begin{theorem}
Let $M$ be a smooth closed oriented $(8n+4)$-manifold and let $F\subset M$ be its characteristic submanifold. Then:
\begin{itemize}
\item
$I(F)\in \Z/8$ depends only on $F$, but not on the auxilliary $\Spin$-structure on $M\sm F$.
\item
$\s(M)=\s(F\circ F)+ 2I(F)\mod 16.$
\end{itemize}
\end{theorem}

\subsection{Realization of quadratic forms}
Rokhlin \cite{R4} pointed out that his congruence answers in negative to 
question \textit {whether every integral unimodular quadratic form
can be realized as an intersection form of a smooth simply connected closed 4-manifold}:
the simplest non-realizable example is the form $E_8$.
And only thirty-years later, further restrictions were obtained
due to the
breakthrough results of S.~Donaldson \cite{Donald83}, who proved
in particular that no positive definite unimodular form other than the standard diagonalizable one
can be realized in such a way.

\subsection{Obstructions to differentiability of manifolds in dimensions $\ge4$}
In Rokhlin's Signature Divisibility Theorem, smoothness of $M$ is crucial, as was shown by M.~Freedman \cite{Freed}:
any unimodular quadratic form is realized by a topological simply-connected 4-manifold.
For instance, Friedman's manifold $M_{E_8}$ realizing form $E_8$  is not smoothable.

Morevover, non-smoothability of $M_{E_8}$ is stable in the sense that
$M_{E_8}\times\R^n$ is also not smoothable
for any $n>0$ (here, the results of Kirby-Siebenmann \cite{KS} are essentially used).

\subsection{Rokhlin's invariant and Rokhlin's homomomorphism}
If $\S$ is an oriented $\Z$-homology 3-sphere,
then, by Theorem \ref{spin-zero}, $\S=\partial W$ for
some smooth oriented 4-manifold $W$ with $w_2(W)=0$.
Then the signature of  $W$ is divisible by 8, and Theorem \ref{divisibility} implies that
the modulo 2 residue of $\sigma(W)/8$, denoted $\mu(\S)\in\Z/2$ and often called the {\it Rokhlin invariant} of $\S$,
is independent of the choice of $W$. It defines the
{\it Rokhlin homomorphism} $\mu : \Theta_3\to \Z/2$ from
the group $\Theta_3$ of
$\Z$-homology 3-sphere classes (where $\S\sim0$ if $\S$ bounds a $\Z$-homology 4-disc).

In 1985 A.~Casson discovered an invariant $\lambda(\S)\in\Z$ of $\Z$-homology spheres
such that $\lambda(\S)=\mu(\S)\mod 2$ and used
it to deduce from Rokhlin's Theorem \ref{divisibility} that  Freedman's manifold $M_{E_8}$ cannot be triangulated
(see Akbulut--McCarthy exposition \cite{AM} for details).

Finally, C.~Manolescu \cite{M}
constructed some integer liftings $\Theta_3\to \Z$ of Rokhlin's $\mu$-homomorphism
which allowed him
to disprove the conjecture on triangulability of all higher-dimensional manifolds
(the idea of proving it using the lifting of $\mu$ goes back to
the foundational works of Casson, Galewski-Stern and Matumoto in 1980s).

\subsection{Exotic Diff structures on a 4-manifold}
A version of Rokhlin's theorem \ref{refined} was used by
S.~Cappell and J.~Shaneson \cite{CS} in their construction of
an exotic projective 4-space (homeomorphic, but not diffeomorphic to $\Rp4$).
Another example
is an exotic $S^3\times \R$ of
M.~Freedman \cite{Fr_fake, Freed}.
It is not diffeomorphic to $S^3\times \R$, since it contains
a Poincar\'e sphere, $P$, which cannot be smoothly embedded
in $S^4\supset S^3\times \R$ because of its Rokhlin invariant $\mu(P)=1$.

%% file: frame-rokhlin.bbl
\begin{thebibliography}{AAAAAA}

\bibitem[GM]{Book}
L.~Guillou, A.~Marin.
\textit{A la recherche de la topologie perdue.}
Progr. Math., 62, Birkhäuser Boston, Boston, MA, 1986.

\bibitem[KM1]{KM2}
M.A.~Kervaire, J.W.~Milnor.
\textit{Bernoulli numbers, homotopy groups, and a theorem of Rohlin.}
1960 Proc. Internat. Congress Math. 1958 pp. 454–458 Cambridge Univ. Press, New York.

 \bibitem[P1]{P1}
 L.S.~Pontryagin.
\textit{Homotopy classification of the mappings of an (n+2)-dimensional sphere on an n-dimensional one.}
Doklady Akad. Nauk SSSR (N.S.)  \textbf{70} (1950), 957 - 959.


 \bibitem[P2]{P2}
L.S.~Pontryagin.
\textit{Characteristic cycles on differentiable manifolds.}
Mat. Sbornik (N. S.)  \textbf{21(63)} (1947), 233 - 284.

 \bibitem[R1]{R1}
 V.A.~Rokhlin.
\textit{On a mapping of an $(n+3)$-dimensional sphere to an $n$-dimensional one.}
Doklady Akad. Nauk SSSR (N.S.)  \textbf{80} (1951), 541 - 544.

 \bibitem[R2]{R2}
 V.A.~Rokhlin.
\textit{Classification of mappings of an $(n+3)$-dimensional sphere to an $n$-dimensional one.}
Doklady Akad. Nauk SSSR (N.S.) \textbf{81} (1951), 19 - 22.

\bibitem[R3]{R3}
 V.A.~Rokhlin.
\textit{A three-dimensional manifold is the boundary of a four-dimensional one.}
Doklady Akad. Nauk SSSR (N.S.) \textbf{81} (1951), 355 - 357.

\bibitem[R4]{R4}
 V.A.~Rokhlin.
\textit{New results in the theory of four-dimensional manifolds.}
Doklady Akad. Nauk SSSR (N.S.) \textbf{84} (1952), 221 - 224.

\bibitem[R5]{R5}
 V.A.~Rokhlin.
\textit{Intrinsic definition of Pontryagin's characteristic cycles.}
Doklady Akad. Nauk SSSR (N.S.) \textbf{84} (1952), 449--452.

\bibitem[R6]{Kolomna}
 V.A.~Rokhlin.
 \textit{Relations between characteristic classes of four-dimensional manifolds.}
 Kolomen. Ped. Inst. Uch. Zap. Ser. Fiz.-Mat. 
\textbf{2} (1958), no. 1, 3 --17.

\bibitem[Se]{SerreNote}
J.A.~Serre.
\textit{Sur les groupes d'Eilenberg-MacLane.}
C. R. Acad. Sci. Paris \textbf{234} (1952), 1243--1245.



\bibitem[Th]{ThomPaper}
R.~ Thom.
\textit{Quelques propri\'et\'es globales des vari\'et\'es diff\'erentiables.}
Comment. Math. Helv. \textbf{28} (1954), 17 – 86.

\bibitem[Wh]{Wh}
H.~Whitney.
\textit{The Singularities of a Smooth $n$-Manifold in $(2n - 1)$-Space.}
Ann. of Math. \textbf{45} (1944), 247--293.

\bibitem[W]{W}
G.W.~Whitehead.
\textit{On the homotopy groups of spheres and rotation groups.}
Ann. of Math.  \textbf{43} (1942), 634--640.


\bibitem[Wu]{Wu}
W.Wu.
\textit{Les $i$-carr\'es dans une vari\'et\'e grassmannienne.}
C. R. Acad. Sci Paris, \textbf{230} (1950), 918--920.

\end{thebibliography}

\begin{thebibliography}{AAAAAA}

\bibitem[AM]{AM}
S.~Akbulut, J.D.~McCarthy.
\textit{Casson's Invariant for Oriented Homology 3-spheres.}
Princeton University Press 1990.


\bibitem[AH]{AH}
M.F.~Atiyah, F.~Hirzebruch.
\textit{Riemann-Roch theorems for differentiable manifolds.}
 Bull. Amer. Math. Soc., \textbf{65} (1959), 276 -- 281.

\bibitem[AS]{AS}
M.F.~Atiyah,  I.M.~Singer.
\textit{The index of elliptic operators on compact manifolds.}
 Bull. Amer. Math. Soc. \textbf{69} (1963), 422–433.

\bibitem[CS]{CS}
S.E.~Cappell, J.L.~Shaneson.
\textit{Some New Four-Manifolds.}
Ann. of Math. \textbf{104} (1976), 61 -- 72.

\bibitem[D1]{Donald83}
S.~Donaldson.
\textit{Self-dual connections and the topology of smooth 4-manifolds.}
Bull. Amer. Math. Soc. (N.S.), \textbf{8} (1983), no. 1, 81 -- 83.


\bibitem[D2]{Donald86}
S.~Donaldson.
\textit{Connections, cohomology and the intersection forms of 4-manifolds.}
J.  Differential  Geom.  \textbf{24}  (1986), 275 -- 341.

\bibitem[Fi]{Finashin}
S.~Finashin.
\textit{A $Pin^-$-cobordism invariant and a generalization of the Rokhlin signature congruence.}
Algebra i Analiz \textbf{2} (1990), no. 4, 242--250.

\bibitem[Fr1]{Fr_fake}
M.H.~Freedman.
\textit{ A fake $S^3\times \R$.}
Ann. of Math. \textbf{ 110} (1979),  177 -- 201.

\bibitem[Fr2]{Freed}
M.H.~Freedman.
\textit{ The topology of four-dimensional manifolds.}
 J. Differential Geometry \textbf{17} (1982), no. 3, 357--453.

\bibitem[FK]{FK}
M.~Freedman, R.~Kirby.
\textit{A geometric proof of Rochlin's theorem.}
in  Algebraic and geometric topology (Proc. Sympos. Pure Math., Stanford Univ., Stanford, Calif., 1976), Part 2, pp. 85–97.

\bibitem[GM-2]{GM}
 L.~Guillou, A.~Marin.
 \textit{ Une extension d'un théorème de Rohlin sur la signature.}
  C. R. Acad. Sci. Paris Sér. A-B
  \textbf{285} (1977), no. 3, A95–A98.



\bibitem[KM2]{KM}
M.A.~Kervaire, J.W.~Milnor.
\textit{On 2-spheres in 4-manifolds.}
Proc. Nat. Acad. Sci. U.S.A. \textbf{47} (1961), 1651 -- 1657.



\bibitem[K]{KirbyLN}
R.~Kirby.
\textit{The topology of 4-manifolds.}
 Lecture Notes in Mathematics, \textbf{1374}. Springer-Verlag, Berlin, 1989.

\bibitem[KS]{KS}
R.~Kirby, L.C.~Siebenmann.
\textit{Foundational essays on topological manifolds, smoothings, and triangulations. With notes by John Milnor and Michael Atiyah.}
 Annals of Mathematics Studies, No. 88. Princeton University Press, Princeton, N.J.; University of Tokyo Press, Tokyo, 1977.

\bibitem[M]{M}
C.~Manolescu.
\textit{ Pin(2)-equivariant Seiberg-Witten Floer homology and the triangulation conjecture.}
J. Amer. Math. Soc. \textbf{29} (2016), no. 1, 147 -- 176.

\bibitem[O1]{O1}
S.~Ochanine.
\textit{The signature of SU-manifolds.}
Mat. Zametki \textbf{13} (1973), 97 --102.

\bibitem[O2]{Ochanine}
S.~Ochanine.
\textit{Signature modulo 16, unvariants de Kervaire generalises et nombres caracteristiques dans la K-theorie reele}
Memoire de la Soc.Math.de France, Nouv.Ser. \textbf{5} (1981).

\bibitem[P]{Pont}
L.S.~Pontryagin.
\textit{Smooth manifolds and their applications in homotopy theory.}
American Mathematical Society Translations, Ser. 2, Vol. 11 (1959),  pp. 1–114.


\bibitem[R7]{R6}
 V.A.~Rokhlin.
\textit{Proof of a conjecture of Gudkov.}
Funkcional. Anal. i Priložen. \textbf{6} (1972), no. 2, 62 -- 64.

 \bibitem[Sco]{Scorpan}
 A.~Scorpan.
 \textit{The wild world of 4-manifolds.}
  American Mathematical Society, Providence, RI, 2005.


%
%
%



\end{thebibliography}
